\documentclass[11pt]{article}
\usepackage[T1]{fontenc}
\usepackage[english]{babel}
\usepackage{verbatim}
\usepackage{amsmath}
\usepackage{amssymb}
\usepackage{xcolor}

\newcommand{\ad}{\mbox{ad}}
\newcommand{\rank}{\mbox{rank}\ }

\newcommand{\demo}{\noindent\textit{Proof. }}

\newcommand{\Z}{\mathcal{Z}}

\newcommand{\C}{\mathbb{C}}
\newcommand{\R}{\mathbb{R}}

\newcommand{\g}{\mathfrak{g}}

\newcommand{\Li}{\mathcal{L}}

\newcommand{\X}{\mathcal{X}}

\newcommand{\ddx}[1]{\frac{\partial}{\partial#1}}
\newcommand{\ds}{\displaystyle}
\newcommand{\pa}{\partial}

\newcommand{\transv}{\mathrel{\text{\tpitchfork}}}
\makeatletter
\newcommand{\tpitchfork}{%
	\vbox{
		\baselineskip\z@skip
		\lineskip-.52ex
		\lineskiplimit\maxdimen
		\m@th
		\ialign{##\crcr\hidewidth\smash{$-$}\hidewidth\crcr$\pitchfork$\crcr}
	}%
}
\makeatother

\begin{document}

\title{General Nilpotent and Solvable Approximations\\ of Almost Riemannian Structures}
\author{Philippe Jouan  and Ronald Manriquez}
\maketitle
\begin{abstract}
The nilpotent or solvable approximation of an ARS is always a linear ARS on a  Lie group or a homogeneous space.
\end{abstract}

\section{Introduction}

We consider here an Almost-Riemannian Structure (ARS in short) defined in a neighborhood $V$ of $p=0$ in $\R^n$ by a set $\{X_1,X_2,\dots,X_n\}$ of n vector fields, considered as an orthonormal frame, which defines the metric. As usual we denote by $\Z$ the singular set, that is the set of points $x\in V$ where the rank of $X_1(x),X_2(x),\dots,X_n(x)$ is not full. The following assumptions guarantee that the ARS is well-defined:
\begin{enumerate}
	\item the set of vector fields $\{X_1,X_2,\dots,X_n\}$ satisfies the rank condition on $V$;
	\item the singular set is not empty but with empty interior.
\end{enumerate}
This implies in particular that $\{X_1,X_2,\dots,X_n\}$ are linearly independant as vector fields.


\subsection{Nilpotent and solvable approximations}\label{NSapprox}

We assume that the point $p=0$ belongs to $\Z$ and that the canonical coordinates of $\R^n$ are privileged with respect to this point. The nilpotent approximation of $X_i$ is denoted by $\widehat{X}_i$.

It is well-known that the  nilpotent approximation $\widehat{X}_1,\ \dots ,\ \widehat{X}_n$ satisfies the rank condition as soon as it is satisfied by $X_1,\dots,X_n$. However it may happen that the vector fields  $\widehat{X}_1,\ \dots ,\ \widehat{X}_n$ fail to be linearly independant, because some of them globally vanish or become a linear combination of the others.

The purpose being here to get an approximation of the original structure that is itself an ARS we proceed as follows:
\begin{enumerate}
	\item Let $k$ be the rank of $\widehat{X}_1,\ \dots ,\ \widehat{X}_n$ at $p=0$. We select a set of $k$ vector fields $\widehat{X}_1,\ \dots ,\ \widehat{X}_k$ linearly independant at $0$.
	\item The remaining vector fields $\widehat{X}_{k+1},\ \dots ,\ \widehat{X}_n$ are assumed to vanish at $0$. If not they can be modified by linear combinations of $\widehat{X}_1,\ \dots ,\ \widehat{X}_k$.
	
	Among them we select a maximal set of linearly independant (as vector fields) elements, say $\widehat{X}_{k+1},\ \dots ,\ \widehat{X}_m$.
	\item If $m<n$ the remaining vector fields either vanish or are linear combinations of $\widehat{X}_{k+1},\ \dots ,\ \widehat{X}_n$. They are replaced by their homogeneous component of degree $0$, denoted by $\widetilde{X}_i$.
\end{enumerate}

\textbf{Remarks}.
\begin{enumerate}
	\item \textit{It may happen that the approximating vector fields obtained by this procedure are not linearly independant. Then the approximating structure is no longer an ARS, it is a sub-Riemannian structure, the rank of which can be constant or not.
	}
	\textit{It seems difficult, if not impossible, to go one step further by considering homogeneous approximations of degree $s>0$ because, as explained in the sequel, two important properties could be lost. First a homogeneous vector field of degree $s>0$ need not be complete. Second the Lie algebra generated by the approximating vector fields would not be finite dimensional in general.  These two drawbacks are related, see [Palais] or [Jouan].}
	
\textit{	We will not consider these degenerate cases in what follows.}
	
	
	\item \textit{It is possible that one, or more, of the $X_i$ has been replaced by a linear combination of the previous ones.  This does not modify neither the dynamics
	$$
	\dot{x}=\sum_{i=1}^{n}u_i X_i,
	$$
	nor the privileged coordinates.}
	
	\textit{However it should be noticed that these substitutions modify the metric since after linear combinations the vector fields are no longer orthonormal but this drawback will be easily solved in the sequel.}
\end{enumerate}


The procedure described above leads to the following set of approximating vector fields:
$$
\widehat{X}_1,\ \dots ,\ \widehat{X}_k,\ \widehat{X}_{k+1},\ \dots,\ \widehat{X}_m,\ \widetilde{X}_{m+1},\ \dots,\ \widetilde{X}_n,
$$
where
\begin{itemize}
	\item $\widehat{X}_i(0)\neq 0$ for $i=1,\dots k$;
	\item $\widehat{X}_i\neq 0$ but $\widehat{X}_i(0)= 0$ for $i=k+1,\dots m$;
	\item $\widehat{X}_i=0$ for $i=m+1,\dots n$.
\end{itemize}

\textbf{In what follows it will always be assumed that these vector fields are independant}. As they are polynomial the singular set $\Z$ where their rank is not full has an empty interior. But $Z$ cannot be empty because at least one of the vector fields vanish at $p=0$. Consequently $\Z$ is a proper with empty interior  subset of $\R^n$ and the set of approximating vector fields defines an ARS.





\subsection{Completness} In view of the next sections it is very important to notice that all involved  vector fields are complete, because of their triangular form.

This fact is well-known for the $\widehat{X}_i$ (see [Bellaïche] or [Jean]). 

It is as well true for the order one homogeneous fields $\widetilde{X}_i$. Indeed such a field writes:
$$
\widetilde{X}_i=\sum_{j=1}^{n}P_{ij}(x)\ddx{x_j}
$$
where $P_{ij}$ is polynomial of nonholonomic degree $w_j$. In particular  $P_{ij}$ depends linearly on the coordinates $(x_{j1},\dots,x_{js})$ of weight $w_j$, is polynomial in the coordinates of weight smaller than $w_j$, and does not depend on the coordinates of weight greater than $w_j$. The associated differential equation is consequently triangular: the coordinates of weight $1$ satisfy a linear homogeneous equation, the coordinates of weight $2$ satisfy a linear equation with a second member that depends on the the coordinates of weight $1$ only, and so on. All solutions are therefore defined on $\R$ and $\widetilde{X}_i$ is complete.

\vskip 0.1cm

This important property does not hold for homogeneous vector fields of positive degree, for example the first coordinate such a vector field could be $\ds x_1^2\ddx{x_1}$, but $\dot{x}_1=x_1^2$ is not complete.


\subsection{The generated Lie algebra}

\newtheorem{FiniteLieAlgebra}{Proposition}
\begin{FiniteLieAlgebra}
	The Lie algebra $\mathcal{L}$ generated by $$
	\widehat{X}_1,\ \dots ,\ \widehat{X}_k,\ \widehat{X}_{k+1},\ \dots,\ \widehat{X}_m,\ \widetilde{X}_{m+1},\ \dots,\ \widetilde{X}_n,
	$$ is finite dimensional.
\end{FiniteLieAlgebra}
\demo
The nonholonomic order of these vector fields is $0$ or $-1$, and the nonholonomic order of their brackets is in the range $-r,\dots,0$ where $r$ is the nonholonomic degree of the set of vector fields. Consequently their components are polynomials of degree less than or equal to $r$. The Lie algebra $\mathcal{L}$ is thus a subspace of a finite dimensional vector space of polynomials.

\hfill $\blacksquare$


\section{The nilpotent case}

It is the case where, after possible linear combinations, the vector fields $\widehat{X}_1,\dots, \widehat{X}_n$ are linearly independant and the vectors $\widehat{X}_1(0),\dots, \widehat{X}_k(0)$ are independant in $\R^n$. In particular no vector field $\widehat{X}_i$ vanishes, and $m=n$.

For $j=k+1,\dots, n$ let $D_j$ stand for $\ad(\widehat{X_j})$ and for any multi-index $J=(j_1,...,j_s)$ let $D_J=D_{j_s}\circ\dots\circ D_{j_1}$ (here $k+1\leq j_u\leq n$ and $s\geq0$). Let
$$
\mathcal{D}=\mbox{Span}\{D_J(\widehat{X}_i)/ \ i=1,\dots p;\ J\ \mbox{ as above} \}.
$$

\newtheorem{Der}{Lemma}
\begin{Der}\label{Der}
The Lie algebra $\g$ generated by $\mathcal{D}$ is $D_j$-invariant for $j=k+1,\dots,n$. 
\end{Der}
\demo
Let $D_{J_1}(\widehat{X}_{i_1})$ and $D_{J_2}(\widehat{X}_{i_2})$ in $\mathcal{D}$. Then
$$
D_j[D_{J_1}(\widehat{X}_{i_1}),D_{J_2}(\widehat{X}_{i_2})]=[D_j\circ  D_{J_1}(\widehat{X}_{i_1}),D_{J_2}(\widehat{X}_{i_2})]+[D_{J_1}(\widehat{X}_{i_1}),D_j\circ D_{J_2}(\widehat{X}_{i_2})]
$$
belongs to $\g$.

\hfill $\blacksquare$

\vskip 0.2cm

Let $\Li$ stand for the Lie algebra generated by $\widehat{X}_1,\ \dots ,\ \widehat{X}_n$. It is a well-known fact that this Lie algebra is nilpotent and finite dimensional (see [Bellaïche] or [Jean]).

\newtheorem{NilpotentCase}{Theorem}
\begin{NilpotentCase} \label{Th1}
	\begin{enumerate}
		\item The ideal generated in $\Li$ by $\widehat{X}_1,\ \dots ,\ \widehat{X}_k$ is $\g$. It is a nilpotent Lie algebra.
		\item The vector fields $\widehat{X}_{k+1},\ \dots ,\ \widehat{X}_n$ do not belong to $\g$ and act on $\g$ as derivations.
		\item The rank at $p=0$ of the elements of $\g$ is full. 
	\end{enumerate}
\end{NilpotentCase}
\demo\begin{enumerate}
	\item Since $\g$ is $D_j$-invariant for $j>k$ it is clear that it is an ideal of $\Li$ that contains $\widehat{X}_1,\ \dots ,\ \widehat{X}_k$. The ideal generated by these vector fields should contain all the $D_J(X_i)$ hence is equal to $\g$.
	\item Let us assume that $\widehat{X}_j$ belongs to $\g$ for some $j>k$. Because of the rule about the nonholonomic order of brackets of homogeneous vector fields, all the elements of $\g$ of order $-1$ are linear combinations of $\widehat{X}_1,\ \dots ,\ \widehat{X}_k$. The vector field $\widehat{X}_j$ is homogeneous of order $-1$ and can consequently be written as:
	$$
	\widehat{X}_j=\sum_{i=1}^{k}\lambda_i \widehat{X}_i.
	$$
	But $\widehat{X}_j$ vanishes at $0$ and the vectors $\widehat{X}_i(0)$ are independent by assumption, so that the $\lambda_i$'s are all equal to $0$, a contradiction.
	\item Let $Y\in\Li\setminus\g$. It can be obtained only as brackets of $\widehat{X}_{k+1},\ \dots ,\ \widehat{X}_n$, hence $Y(0)=0$ because all these fields vanish at $0$.
	
	On the other hand the Lie algebra $\Li$ satisfies the rank condition, hence $\g$ should satisfy the rank condition at $0$.
\end{enumerate}

\hfill $\blacksquare$

\vskip 0.1cm

After this analysis at the algebra level we can turn our attention to the Lie group level.

Let $G$ be the simply connected Lie group whose Lie algebra is $\g$. Since $\g$ is nilpotent the underlying manifold of $G$ is $\R^N$, $N=\dim(\g)$. The first task is to show that $\R^n$ is a homogeneous space of $G$. This is mainly due to the fact that $\g$ is generated by homogeneous vector fields.

\newtheorem{HomSpace}[Der]{Lemma}
\begin{HomSpace}\label{HomSpace}
	The set $\R^n$ is a homogeneous space of $G$. More accurately if $H$ stands for the connected subgroup of $G$ whose Lie algebra is the set of elements of $\g$ that vanish at $0$, then $\R^n$ is diffeomorphic to the quotient $G/H$.
\end{HomSpace}
\demo

Since the elements of $\g$ are complete vector fields of $\R^n$ the group $G$ acts naturally on $\R^n$ as a group of diffeomorphisms, and it is enough to show that this action is transitive. More accurately $G$ is the set of $\exp(t_sY_s)\dots\exp(t_1Y_1)$ where $Y_i\in\g$ is a vector field on $\R^n$ and $t_i\in\R$, i.e. a group of diffeomorphisms of $\R^n$. 

Let us begin by a simple remark. Let $\ds Y=\sum_{j=1}^{n}y_j(x)\ddx{x_j}$ be an element of $\g$, and assume that $Y$ is homogeneous and that $Y(0)\neq 0$. There exists an index $i$ such that $y_i(0)\neq 0$, and by homogeneity $y_i$ is constant: $\forall x\in\R^n$, $y_i(x)=a_i\neq 0$. By homogeneity again the polynomials $y_j$ are homogeneous of order $w_j-w_i$. Consequently
$$
Y=\sum_{w_j=w_i}a_j\ddx{x_j}+\sum_{w_j>w_i}y_j(x)\ddx{x_j} \quad \mbox{ where }\ \ a_i\neq 0.
$$
Since $\g$ is generated by homogeneous elements and its rank is full at $0$, we can choose in $\g$ a set $n$ homogeneous elements $Y_1,Y_2,\dots, Y_n$ linearly independant at $0$. We can also assume that the order of $Y_1,\dots,Y_{n_1}$ is $-1$, the order of $Y_{n_1+1},\dots,Y_{n_2}$ is $-2$ and so on. Up to linear combinations we can assume that  $Y_1,\dots,Y_{n_1}$ have the following form:
$$
Y_i=\ddx{x_i}+\sum_{w_j>1}y_j(x)\ddx{x_j}.
$$ 
More generally we can assume that if the order of $Y_i$ is $w_i$ then
$$
Y_i=\ddx{x_i}+\sum_{w_j>w_i}y_j(x)\ddx{x_j}.
$$
This way it is clear that the rank of the set of vector fields $Y_1,Y_2,\dots, Y_n$ is full everywhere in $\R^n$. This implies that the action of $G$ on $\R^n$ is transitive

\hfill $\blacksquare$

\vskip 0.2cm

To complete the construction we associate to the derivation $D_j=\ad(\widehat{X}_j)$ of $\g$ a linear vector field $\X_j$ on $G$ for  $j>k$ ($\X_j$ does exist because $G$ is simply connected). It is clear that the projection of $\X_j$ on $\R^n$ is $\widehat{X}_j$ (See \cite{Jouan09} for details). Finally the vector fields $\widehat{X}_1,\ \dots ,\ \widehat{X}_k$ are invariant, and $\widehat{X}_{k+1},\ \dots ,\ \widehat{X}_n$ are linear vector fields on the homogeneous space $\R^n=G/H$.

We can state:

\newtheorem{NilpotentGroup}{Corollary}
\begin{NilpotentGroup} \label{NilpotentGroup}
	The space $\R^n$ is a homogeneous space of the nilpotent Lie group $G$ whose Lie algebra is $\g$.
	
	The vector fields $\widehat{X}_1,\ \dots ,\ \widehat{X}_k$ are projections of invariant vector fields of $G$ and $\widehat{X}_{k+1},\ \dots ,\ \widehat{X}_n$ are projections of linear vector fields of $G$.
	
	Consequently the set $\widehat{X}_1,\ \dots ,\ \widehat{X}_n$ defines a linear ARS on the homogeneous space $\R^n$.
\end{NilpotentGroup}

If during the procedure of Section \ref{NSapprox} one (or more) of the $\widehat{X}_j$, with $j>k$, have been modified by a linear combination of the $\widehat{X}_i$, $i\leq k$, we can now go back to the original approximating vector fields. For this we just have to add to the linear vector field $\widehat{X}_j$ a linear combination of the invariant vector fields $\widehat{X}_1,\dots,\widehat{X}_k$. The resulting vector field is then an affine one, that is the sum of a (nonzero) linear vector field with an invariant one.

\newtheorem{NilpotentGroup2}[NilpotentGroup]{Corollary}
\begin{NilpotentGroup2} \label{NilpotentGroup2}
	The space $\R^n$ is a homogeneous space of the nilpotent Lie group $G$ whose Lie algebra is $\g$.
	
	The vector fields $\widehat{X}_1,\ \dots ,\ \widehat{X}_k$ are projections of invariant vector fields of $G$ and $\widehat{X}_{k+1},\ \dots ,\ \widehat{X}_n$ are projections of linear or affine, not invariant, vector fields of $G$.
	
	Consequently the set $\widehat{X}_1,\ \dots ,\ \widehat{X}_n$ defines a linear ARS on the homogeneous space $\R^n$.
\end{NilpotentGroup2}

\vskip 0.2cm

\noindent \textbf{Example 1}

Let in $\R^3$:
$$
X_1=\pa_x, \quad X_2=x\pa_y, \quad X_3=y^2\pa_z.
$$
The Lie algebra $\Li$ contains $X_1,X_2,X_3$ and
$$
\begin{array}{l}
X_4=[X_1,X_2]=\pa_y,\ X_5=\frac{1}{2}[X_4,X_3]=y\pa_z,\\
X_6=[X_2,X_5]=x\pa_z, \ X_7=[X_1,X_6]=[X_4,X_5]=\pa_z\\
\X_1=\frac{1}{2}[X_2,X_3]=xy\pa_z, \ \X_2=[X_2,\X_1]=x^2\pa_z.
\end{array}
$$
$$
\Delta=\{X_1,X_2,X_3\} \quad \Delta^2=\{X_4,\X_1\} \quad \Delta^3=\{X_5,\X_2\} \quad \Delta^4=\{X_6\} \quad \Delta^5=\{X_7\}
$$
The canonical coordinates are privileged with weights $(1,2,5)$ and the vector fields $X_1$, $X_2$, and $X_3$ are homogeneous of order $-1$ hence equal to their nilpotent approximations.

The algebra $\g$ is here the ideal of $\Li$ generated by $X_1$ that is
$$
\g=\mbox{Span}\{X_1,X_4,X_5,X_6,X_7 \}.
$$
The vector fields $X_2$ and $X_3$ are linear, as well as the fields $\X_1=\frac{1}{2}[X_2,X_3]$ and $\X_2=[X_2,\X_1]$.

The orders of $\X_1$ and $\X_2$ are respectively $-2$ and $-3$ which shows that the vector fields of order smaller than $-1$ are not necessarily in $\g$.

Notice that the singular locus is here $\Z=\{xy^2=0 \}$.

\hfill $\blacksquare$


\section{The non-nilpotent case}
\subsection{Examples}

In the general case where $m<n$ some difficulties may happen, as shown by the following example.

\vskip 0.2cm

\noindent \textbf{Example 2}

Consider in $\R^4$:
$$
X_1=\pa_x, \quad X_2=\pa_y+x\pa_z, \quad X_3=y\pa_w, \quad X_4=x\pa_z+\frac{1}{2}y^2\pa_w.
$$
The Lie algebra $\Li$ contains also
$$
X_5=[X_1,X_2]=\pa_z,\ X_6=[X_2,X_3]=\pa_w.
$$
The canonical coordinates $(x,y,z,w)$ are privileged with weights $(1,1,2,2)$. The vector fields $X_1$, $X_2$ are homogeneous of order $-1$ and independent at $0$, and the vector field $X_3$ is homogeneous of order $-1$ but vanishes at $0$. The last field $X_4$ is homogeneous of order $0$. Consequently the first three are equal to their nilpotent approximation and $X_4=\widetilde{X}_4$. According to the notations of the introduction we have $k=2$ and $m=3$.

We cannot assert as in the nilpotent case that $X_3$ does not belong to $\g$. Indeed the ideal generated in $\Li$ by $X_1$ and $X_2$ is here
$$
\g=\mbox{Span}\{X_1,X_2,X_3,X_4,X_5 \}
$$
because $X_3=[X_2,X_4]$.

This is a drawback that can happen when $k<m<n$

Notice that the singular locus is here $\Z=\{xy=0 \}$.

\hfill $\blacksquare$

\vskip 0.2cm

In the general case the Lie algebra $\Li$ need not be solvable. Indeed it is a subalgebra of the semi-direct product of $\g$ by its algebra of derivations. But the algebra of derivations of a nilpotent Lie algebra is not solvable in general. For instance the derivations of the Heisenberg algebra is the set of endomorphisms whose matrix in the canonical basis is:
$$
D=\begin{pmatrix}
a & b & 0 \\
c & d & 0 \\
e & f & a+d
\end{pmatrix}
$$
The subalgebra of such derivations that moreover satisfy $e=f=a+d=0$ is equal to $\mathfrak{sl_2}$ hence semisimple.

Example 3 illustrates that phenomenon.

\vskip 0.2cm

\noindent \textbf{Example 3}

Consider in $\R^5$, with coordinates $(x,y,z,w,t)$:
$$
\begin{array}{l}
X_1=\pa_x, \quad  X_2=\pa_y+x\pa_z+w\pa_t, \quad  X_3=\pa_w+x\pa_t, \\  X_4=x\pa_y+\frac{1}{2}x^2\pa_z, \quad  X_5=y\pa_x+\frac{1}{2}y^2\pa_z
\end{array}
$$
The Lie algebra $\Li$ contains also
$$
\begin{array}{l}
X_6=[X_1,X_2]=\pa_z,\quad X_7=[X_1,X_3]=[X_3,X_2]=\pa_t,\\
X_8=[X_3,X_5]=-y\pa_t,\quad \X_9=[X_8,X_4]=x\pa_t,\\
 \X_{10}=[X_4,X_5]=x\pa_x-y\pa_y.
\end{array}
$$
The coordinates $(x,y,z,w,t)$ are privileged with weights $(1,1,2,1,2)$. The vector fields $X_4$, $X_5$ and $\X_{10}$ are homogeneous of order $0$ and does not belong to $\g$. A simple computation shows that they are equivalent to the usual basis of $\mathfrak{sl_2}$. Consequently the algebra $\Li$ is not solvable.

The singular locus is here $\Z=\{??? \}$.

\hfill $\blacksquare$


\subsection{Statements}

We set $D_j=\ad(\widehat{X_j})$ for $j=k+1,\dots, m$ and $D_j=\ad(\widetilde{X_j})$ for $j=m+1,\dots, n$. As well as in the nilpotent case we set $D_J=D_{j_s}\circ\dots\circ D_{j_1}$ for any multi-index $J=(j_1,...,j_s)$ where $s\geq0$ and $k+1\leq j_u\leq n$, and
$$
\mathcal{D}=\mbox{Span}\{D_J(X_i)/ \ i=1,\dots k;\ J\ \mbox{ as above} \}.
$$
The Lie algebra $\g$ generated by $\mathcal{D}$ is again $D_j$-invariant for $j=k+1,\dots, n$, which shows that $\g$ is the ideal generated in $\Li$ by $\widehat{X}_1,\ \dots ,\ \widehat{X}_k$.
 
\newtheorem{NonilpotentCase}[NilpotentCase]{Theorem}
\begin{NonilpotentCase}
	\begin{enumerate}
		\item The ideal generated in $\Li$ by $\widehat{X}_1,\ \dots ,\ \widehat{X}_k$ is $\g$. It is a nilpotent Lie algebra.
		\item The vector fields $\widetilde{X}_{m+1},\ \dots ,\ \widetilde{X}_n$ do not belong to $\g$ and act on $\g$ as derivations.
		\item The vector fields $\widehat{X}_j$, with $k+1\leq j\leq m$ that do not belong to $\g$ act on $\g$ as derivations.
		\item The rank at $p=0$ of the elements of $\g$ is full. 
	\end{enumerate}
\end{NonilpotentCase}
\demo The proofs that $\g$ is the ideal generated in $\Li$ by $\widehat{X}_1,\ \dots ,\ \widehat{X}_k$ and that its rank at $p=0$ is full are identical to the nilpotent case.

The Lie algebra $\g$ is generated by homogeneous vector fields of order at most $-1$. Since the order of a Lie bracket is the sum of the orders of the factors and a vector field of order less than $-r$ vanishes, all brackets of length larger than $r$ vanish, which shows that $\g$ is nilpotent.

The points 2. and 3. are clear.

\hfill $\blacksquare$

\vskip 0.2cm

Because of the phenomenon described in Example 2 we cannot guarantee that the vectors fields $\widehat{X}_{k+1},\ \dots ,\ \widehat{X}_m$ do not belong to $\g$ and we are lead to introduce one more index. Up to a reordering we can assume that $\widehat{X}_{k+1},\ \dots ,\ \widehat{X}_l$ belong to $\g$ and that $\widehat{X}_{l+1},\ \dots ,\ \widehat{X}_m$ do not belong to $\g$, where $k+1\leq l\leq m$.

\newtheorem{NonilpotentGroup}[NilpotentGroup]{Corollary}
\begin{NonilpotentGroup} \label{NonilpotentGroup}
	The space $\R^n$ is a homogeneous space of the nilpotent Lie group $G$ whose Lie algebra is $\g$.
	
	The vector fields $\widehat{X}_1,\ \dots ,\ \widehat{X}_l$ are projections of invariant vector fields of $G$. The vector fields $\widehat{X}_{l+1},\ \dots ,\ \widehat{X}_m$ are projections of linear or affine vector fields of $G$ and $\widetilde{X}_{m+1},\ \dots ,\ \widetilde{X}_n$ are projections of linear ones.
	
	Consequently the set $\widehat{X}_1,\ \dots ,\ \widehat{X}_n$ defines a linear ARS on the homogeneous space $\R^n$.
\end{NonilpotentGroup}

\vskip 0.2cm

\noindent \textbf{Remarks}.

\begin{enumerate}
	\item If $m=n-1$ then the Lie algebra $\Li$ is solvable. It is why we had previously called \textit{solvable} the approximations at order $0$. It is no longer true when $m\leq n-2$, and no general conclusion holds in that case, but we keep the name \textit{solvable approximation}.
	\item 
\end{enumerate}


\section{Genericity}

We consider in this section a distribution $\Delta$ on a $n$-dimensional connected manifold $M$. This distribution is assumed to be locally defined by $n$ vector fields. It is not useful to add that the distribution cannot be locally defined by less that $n$ vector fields because this condition will appear to be generic.

\texttt{\color{red} For the first two theorems we need neither metric, nor normal forms, they deal with distributions only.}


\subsection{Generic distributions}

\newtheorem{Generic1}{Theorem}
\begin{Generic1}
	There exists an open and dense subset of ?? that verify:
	\begin{enumerate}
		\item Let $R$ be the largest integer such that $R^2\leq n$. For $1\leq r\leq R$
		let $\Z_r$ be the set of points where the rank of $\{X_1,X_2,\dots,X_n\}$ is $n-r$. Each $\Z_r$ is a codimension $r^2$ submanifold and the singular locus $\Z$ is the union of these disjoint submanifolds.
		\item The submanifold $\Z_{r+1}$ is included in the closure $\overline{\Z}_r$ of $\Z_r$ for $r=1,\dots, R-1$.
		\item The mapping $x\longmapsto \det(X(x))$ is a submersion at all points $x\in\Z_1$.
		\item The rank of $\Delta+[\Delta,\Delta]$ is full at all points.
	\end{enumerate}
\end{Generic1}

\newtheorem{ParallelPoints}[Generic1]{Theorem}
\begin{ParallelPoints}
	We have	
	\begin{itemize}
		\item $r=1$. There exist isolated points in $\Z_1$ where $T_p\Z_1=\Delta_p$
		\item $r\geq 2$. Let $m(n,r)$ be the largest dimension that $T_p\Z_r+\Delta_p$ may reach, that is $m(n,r)=\min(n,2n-r^2-r)$, and let $s=m(n,r)-\dim(T_p\Z_r+\Delta_p)$. Then
		\begin{enumerate}
			\item The set of points $p\in \Z_r$ where $s=1$ is a submanifold of $\Z_r$ as soon as $\ds n\geq r^2+r-E\left(\frac{r-1}{2}\right)$.
			\item The set of points $p\in \Z_r$ where $s\geq 2$ and $s^2\leq r$  is a submanifold of $\Z_r$ as soon as $\ds r^2+r-E\left(\frac{r-s^2}{s-1}\right)\leq n \leq r^2+r+E\left(\frac{r-s^2}{s+1}\right)$.
			\item The set of points $p\in \Z_r$ where $s\geq 2$ and $s^2>r$ is empty.
		\end{enumerate}
	\end{itemize}
\end{ParallelPoints}


\subsection{Normal forms of generic distributions}

\texttt{\color{red} At this point it is necessary to introduce the almost-Riemannian metric and the normal forms.}


\subsection{Nilpotent and solvable approximations of generic distributions}

\newtheorem{TildeExists}[Generic1]{Theorem}
\begin{TildeExists}
	For a generic distribution holds:
	\begin{enumerate}
		\item[(i)] Let $p$ be a parallel point in $\Z_1$, that is a point where $T_p\Z_1=\Delta_p$. Then $\widehat{X}_n=0$ but $\widetilde{X}_n\neq 0$, in normal form.
		\item[(ii)] At all other points, including the points in $\Z_r;\ \ r\geq 2$, the nilpotent approximation $\widehat{X}_1,\widehat{X}_2,\dots,\widehat{X}_n$ is a set of $n$ linearly independant vector fields.
	\end{enumerate}	
\end{TildeExists}

\section{ExGenericity}

It is clear from the previous sections that many different, complicated structures may happen and that it is of interest to limit ourselves to generic cases, or at least to understand what they are.

\vskip 0.1cm

\textbf{Remark}. An extensive use of the product of coranks theorem (see \cite{GG73}) will be done. Let us recall it.

\textit{Let $\mathcal{M}(n,m)$ be the set of $n\times m$ real matrices. The set $L^r$ of such matrices of corank $r$ is a submanifold  of $\mathcal{M}(n,m)$ of codimension $(n-q+r)(m-q+r)$ where $q=\min\{n,m\}$.}


\subsection{The singular locus}

Let $\{X_1,X_2,\dots,X_n\}$ be a set of $n$ vector fields of class $\C^\infty$ defined in a neighborhood $V$ of $0$ in $\R^n$. We denote by $X$ the $n\times n$ matrix whose $j^{\mbox{th}}$ column is $X_j$. This matrix varies with the point $x\in V$ and $X$ can be viewed as a $\C^\infty$ mapping from $V$ to $\mathcal{M}(n,n)$. According to the Thom's tranversality Theorem (see \cite{GG73}) the set of such mappings that are transversal to $L^r$ is residual in $\C^\infty(V,\mathcal{M}(n,n))$. The codimension of $L^r$ is here $r^2$ and we get:
\begin{itemize}
	\item If $r^2>n$ then $X\transv L^r$ means $X(V)\cap L^r=\emptyset$.
	\item If $r^2\leq n$ then $X\transv L^r$ means that $X^{-1}(L^r)$, if not empty, is a codimension $r^2$ submanifold of $V$.
\end{itemize}

\textbf{Conclusion}. \textit{Let $R$ be the largest integer such that $R^2\leq n$. For $1\leq r\leq R$
let $\Z_r$ be the set of points where the rank of $\{X_1,X_2,\dots,X_n\}$ is $n-r$. It is generically a codimension $r^2$ submanifold of $V$, possibly empty, and the singular locus $\Z$ is the union of the $\Z_r$, hence a union of submanifolds.}

\vskip 0.2cm

\textbf{Remark}. For $n=2,3$ the singular locus is generically a codimension one submanifold equal to $\Z_1$ (it is the statement of \cite{ABS08}, \cite{BCG13}, and \cite{BCGM14})


\subsection{The step of the distribution}

In this section $\Delta_p$ stands for the evalution at $p\in V$ of the distribution $\Delta$ generated by the $X_j$s.
Let us denote by $X_{ij}$ the Lie bracket $[X_i,X_j]$ for $i<j$. There are $\ds \frac{(n-1)n}{2}$ such brackets.

Let $\ds\mathcal{M}=\mathcal{M}(n,n+\frac{(n-1)n}{2})$ and $\Phi$ be the mapping from $V$ to $\mathcal{M}$ whose columns are the $X_j$ and the $X_{ij}$. The codimension of $L^1$ in $\mathcal{M}$ is
$$
(n-n+1)(n+\frac{(n-1)n}{2}-n+1)=\frac{(n-1)n}{2}+1
$$
For $n\geq 3$ this codimension is strictly larger than $n$ and transversality to $L^1$ means nonintersection. The same conclusion holds obviously for $L^r$, $r>1$.

The dimension $2$ is particular: $\Phi^{-1}(L^1)$ is generically a codimension $2$ submanifold of $V$, i.e. there are isolated points where the rank of $\Delta_p+[\Delta_p,\Delta_p]$ is not full (see \cite{ABS08} and \cite{BCG13}).

\textbf{Conclusion}. \textit{Generically, and for $n\geq 3$, holds
	$$
	\rank \left(\Delta_p+[\Delta_p,\Delta_p]\right) =n
	$$
at all points $p$ of $V$.}


\subsection{The tangential points}

\textit{\textbf{Take care}. The proofs in this section are not rigorous enough and must be improved. But the results are probably right...}

\vskip 0.2cm

It is shown in \cite{ABS08} and \cite{BCGM14} that generically, in dimensions $2$ and $3$, the points $p$ belonging to the singular locus but where the distribution is equal to the tangent space to $\Z$ are isolated. Such points are called tangential.

Consider a generic local ARS for $n>3$ and let $p\in \Z_r$. The rank of the distribution at $p$ is by definition $n-r$. On the other hand the dimension of $\Z_r$ is $n-r^2$, strictly less than $n-r$ as soon as $r>1$. We can consequently consider two opposite properties:
\begin{enumerate}
	\item Inclusion, that is $T_p\Z_r\subset \Delta_p$.
	\item Transversality, $T_p\Z_r \transv \Delta_p$.
\end{enumerate}
These two properties are identical for $r=1$ only.

\vskip 0.2cm

Since $\Z_r$ is a dimension $(k=n-r^2)$ submanifold of $V$ we can choose coordinates such that (around $p$):
$$
\Z_r=\{x_{k+1}=\dots=x_n=0 \} \quad \mbox{and }\ \ T_p\Z=\mbox{Span}\{\partial_{x_1},\dots,\partial_{x_k} \}.
$$
We can also select from the distribution $n-r$ vectors fields independent at $p$ and assume without loss of generality that they are $X_1,\dots,X_{n-r}$.

Let $\Phi$ be the mapping from $\Z_r$ to $\mathcal{M}(r^2,n-r)$ defined by
$$
\Phi(q)=\left(\overline{X_1}(q),\dots, \overline{X_{n-r}}(q) \right)
$$
where $\overline{X_j}$ is made of the $r^2$ last coordinates of $X_j$.

The first property, inclusion, holds at $q\in \Z_r$ if the rank of $\Phi(q)$ is $n-r-(n-r^2)=r^2-r$.

Let $s=\min\{r^2-(r-r^2),n-r-(r^2-r) \}=\min\{r,n-r^2 \}$. Then there is inclusion at $q$ if $\Phi(q)\in L^s$. The codimension of $L^s$ in $\mathcal{M}(r^2,n-r)$ is $r(n-r^2)$. For the two conditions $\Phi(q)\in L^s$ and $\Phi$ transversal to $L^s$ to hold it is necessary that $\dim \Z_r\geq \mbox{codim} L^s$, that is $n-r^2\geq r(n-r^2)$. This condition is possible if and only if $n-r^2=0$ or $r=1$. In the first case $\Z_r$ consists of isolated points and the inclusion is natural, in the second one the codimension of $L^s$ is $n-1$ and the set of tangential points is a codimension $n-1$ submanifold of $\Z_1$ hence made of isolated points. This result is coherent with the ones of \cite{ABS08}, \cite{BCG13}, and \cite{BCGM14}.

\vskip 0.2cm

The second property, transversality to $\Z_r$, is not always possible. Indeed it implies $\dim(\Z_r)+\dim\Delta_p\geq n$, that is $n-r^2+n-r\geq n$ or $n\geq r^2+r$. 

\textbf{And is it useful?}


\section{Generic normal forms}

Despite the fact that our purpose is local it will be more convenient here to define an ARS by a smooth distribution $\Delta$ endowed with an inner product, smoothly varying with the point.

All along this section the distribution will be assumed generic, that is:
\begin{enumerate}
	\item The singular locus $\Z$ is a union of submanifolds $\Z_r$ of codimension $r^2$ where the rank of $\Delta$ is $n-r$.
	\item $	\rank \left(\Delta_p+[\Delta_p,\Delta_p]\right) =n$ at all points.
	\item The tangency (or inclusion) points are in $\Z_1$ only and isolated. 
	\end{enumerate}

First we can follow the lines of $\cite{ABS08}$ (also used in \cite{BCG13} and \cite{BCGM14}).

Let $W$ be a codimension 1 submanifold transversal to the distribution. We can define a coordinate system $y=(x_2,\dots,x_n)$ in $W$, and choose an orientation transversal to $W$. Let $\gamma_y$ be the family of normal geodesics parametrized by arclength, transversal to $W$ at $y$, and positively oriented. The mapping $(x_1,y)\mapsto \gamma_y(x_1)$ is a local diffeomorphism and the geodesics $x_1\mapsto \gamma_y(x_1)$ realize the minimal distance between $W=\{ x_1=0\}$ and the surfaces $\{x_1=c\}$ for $c$ small enough. The transversality conditions of the PMP are consequently satisfied along all these surfaces: if $\lambda(x_1)$ is a covector associated one of these geodesics then the tangent space to $\{x_1=c\}$ at $\gamma_y(x_1)$ is $\ker(\lambda(x_1))$

Now let $X_1=\partial_{x_1}$ be the vector field defined in the following way:  $X_1(q)=\ds \frac{d}{d x_1}\gamma_y(x_1)$ at the point $q=\gamma_y(x_1)$. It is a unitary vector field belonging to $\Delta$. Let $\{X_1,X_2,\dots,X_n\}$ be any orthonormal frame of $\Delta$. Since the controls of the geodesics $\gamma_y$ are $(1,0\dots,0)$ and $u_j=\left\langle \lambda, X_j\right\rangle $ the vector fields $X_2,\dots,X_n$ are tangent to the surfaces $\{x_1=c\}$.. Consequently the vector fields have the following form:
$$
X_1=\begin{pmatrix} 1 \\ 0 \\  . \\ . \\ 0 \end{pmatrix},\quad
X_j=\begin{pmatrix} 0 \\ b_{2,j} \\ .\\  . \\ b_{n,j} \end{pmatrix}\ \ \mbox{ for } 1<j<n \ \ 
$$
for any choice of the coordinates in $W$ and any choice of $X_2,\dots,X_n$, under the condition that they provide an orthonormal frame related to the subRiemannian metric.

\vskip 0.2cm

\noindent \textbf{Important remark}
\textit{Notice the vectors of the distribution that are orthogonal to $X_1$ \textbf{for the subRiemannian metric } are orthogonal to $X_1$ for the canonical inner product of $\R^n$ \textbf{for the chosen coordinates}. This is due to the transversality conditions.}

\vskip 0.2cm

Let us consider the case where $p=0$ belongs to $\Z_1$ and is not a tangent point. We can then choose $W=\Z_1$. We obtain a parameter $x_1$ of the geodesics from $\Z_1=W=\{x_1=0\}$ to $\{x_1=c\}$ for $c$ small enough.

Certainly it is possible to show that the vector fields have the following form in suitable coordinates:
$$
X_1=\begin{pmatrix} 1 \\ 0 \\ . \\  . \\ . \\ . \\ 0 \end{pmatrix},
X_j=\begin{pmatrix} 0 \\ . \\ 0 \\  1+a_{jj} \\ b_{j+1,j} \\ . \\ b_{n,j} \end{pmatrix}\ \ \mbox{ for } 1<j<n \ \ \mbox{ and }
X_n=\begin{pmatrix} 0 \\ . \\ . \\  . \\ . \\ 0 \\ a_{n,n} \end{pmatrix}.
$$
where $a_{jj}(0)=0$ for $j>1$ and $b_{i,j}(0)=0$ for $1<j<i\leq n$.
Moreover $a_{n,n}=0\Longleftrightarrow x_1=0$.

I have not really proved that but I think it is true. It would  of great interest to prove it and then to analyze the other generic cases: the tangent case and then the points where the rank of the distribution is strictly less than $n-1$.

To finish we need to know in which generic cases it may happen that some $\widehat{X}_j$ vanishes and must be replaced by $\widetilde{X}_j$. It is even our main purpose.

\vskip 0.3cm


\end{document}